\documentclass[11pt]{article}
\usepackage{amssymb,amsmath,latexsym, verbatim}
\allowdisplaybreaks

\begin{document}

\begin{doublespace}

\newtheorem{thm}{Theorem}[section]
\newtheorem{lemma}[thm]{Lemma}
\newtheorem{cond}[thm]{Condition}
\newtheorem{defn}[thm]{Definition}
\newtheorem{prop}[thm]{Proposition}
\newtheorem{corollary}[thm]{Corollary}
\newtheorem{remark}[thm]{Remark}
\newtheorem{example}[thm]{Example}
\newtheorem{conj}[thm]{Conjecture}
\numberwithin{equation}{section}
\def\ee{\varepsilon}
\def\qed{{\hfill $\Box$ \bigskip}}
\def\NN{{\cal N}}
\def\AA{{\cal A}}
\def\MM{{\cal M}}
\def\BB{{\cal B}}
\def\CC{{\cal C}}
\def\LL{{\cal L}}
\def\DD{{\cal D}}
\def\FF{{\cal F}}
\def\EE{{\cal E}}
\def\QQ{{\cal Q}}
\def\RR{{\mathbb R}}
\def\R{{\mathbb R}}
\def\L{{\bf L}}
\def\K{{\bf K}}
\def\S{{\bf S}}
\def\A{{\bf A}}
\def\E{{\mathbb E}}
\def\F{{\bf F}}
\def\P{{\mathbb P}}
\def\N{{\mathbb N}}
\def\eps{\varepsilon}
\def\wh{\widehat}
\def\wt{\widetilde}
\def\pf{\noindent{\bf Proof.} }
\def\beq{\begin{equation}}
\def\eeq{\end{equation}}
\def\lam{\lambda}
\def\H{\mathcal{H}}
\def\nn{\nonumber}
\def\L{\mathcal{L}}

\newcommand{\Per}{\mathrm{Per}}
\newcommand{\norm}[1]{{\lVert #1 \rVert}}

\title{\Large \bf  Small time asymptotics of spectral heat contents for subordinate killed Brownian motions related to 
isotropic $\alpha$-stable processes}
\author{ \bf  Hyunchul Park\hspace{1mm} \hspace{1mm}\hspace{1mm} and \hspace{1mm}\hspace{1mm}
Renming Song\thanks{Research supported in part by a grant from the Simons
Foundation (\#429343, Renming Song).}\hspace{1mm} }

\date{ \today}
\maketitle

\begin{abstract}
In this paper we study the small time asymptotic behavior of the spectral heat content 
$\widetilde{Q}_D^{(\alpha)}(t)$ of an arbitrary bounded $C^{1,1}$ domain $D$  with respect to the \textit{subordinate killed Brownian motion} in $D$ via an $(\alpha/2)$-stable subordinator. 
For all  $\alpha\in (0,2)$, we establish a two-term small time expansion for $\widetilde{Q}_D^{(\alpha)}(t)$ in all dimensions. 
When $\alpha\in (1,2)$ and $d\geq 2$, we establish a three-term small time expansion for
$\widetilde{Q}_D^{(\alpha)}(t)$.
\end{abstract}

\medskip


\medskip


\bigskip
\section{Introduction}\label{introduction}

\quad 
Let $X$ be a Markov process in $\R^d$. For any open set $D\subset \R^d$, the heat content
of $D$ with respect to $X$ is defined to be
$$
H^X_D(t):=\int_D\P_x(X_t\in D)dx,
$$
and the spectral heat content of $D$ with respect to $X$ is defined to be
$$
Q^X_D(t):=\int_D\P_x(\tau^X_D>t)dx,
$$
where $\tau^X_D$ is the first time the process $X$ exits $D$. 
The spectral heat content of $D$ with respect to $X$ can be regarded as the heat content of $D$ with respect to the killed
process $X^D$. When $X$ is an isotropic $\alpha$-stable process, $\alpha\in (0, 2]$, in
$\R^d$, we will write  $Q^{(\alpha)}_D(t)$ for $Q^X_D(t)$. In particular, $Q^{(2)}_D(t)$
stands for the spectral heat content of $D$ with respect to Brownian motion.

The heat content with respect to  L\'evy processes, especially Brownian motions, has been
studied extensively, see, for instance, \cite{Valverde2}, \cite{V1}, \cite{vanDenBerg1_POT}, \cite{vanDenBerg1}, \cite{vanDenBerg2}, and \cite{cg}.  
The spectral heat content $Q^{(2)}_D(t)$ with respect to Brownian
motion has also been studied a lot (see \cite{V2}--\cite{vanDenBerg3}, \cite{BD}, and \cite{BG}).  
In \cite{BD} a two-term small time expansion for $Q^{(2)}_D(t)$ was established for bounded $C^{1,1}$ domains 
and in \cite{BG} a three-term small time expansion for $Q_{D}^{(2)}(t)$ was obtained for bounded domains with $C^{3}$ boundary.
In \cite{Savo} a recursive formula of the complete asymptotic series of the spectral heat content in a Riemannian manifold with smooth boundary was investigated. The study of the
small time asymptotic behavior of the  spectral heat content with respect to other
L\'evy processes is more recent.  
Upper and lower bounds for $Q_D^{(\alpha)}(t)$, $\alpha\in (0, 2)$, were established in 
\cite{V1}, while explicit expressions for the second term in the asymptotic behavior of $Q_{D}^{(\alpha)}(t)$, $\alpha\in (0,2)$, in dimension 1 for bounded open intervals were obtained in \cite{V2}.
In the recent paper \cite{GPS}, 
the results of \cite{V1} and \cite{V2} were generalized in several directions.

An isotropic  $\alpha$-stable process $X^{(\alpha)}$ can be obtained from a Brownian motion $W$
via an independent  $(\alpha/2)$-stable subordinator $S^{(\alpha/2)}$, 
that is, $X^{(\alpha)}_t=W_{S^{(\alpha/2)}_t}$. Thus
an isotropic  $\alpha$-stable process is a subordinate Brownian motion.
Hence the spectral heat content $Q^{(\alpha)}_D(t)$
is the heat content with respect to the \textit{killed subordinate Brownian motion} $X^{(\alpha), D}$, which 
can be obtained from the Brownian motion $W$ by subordinating with the independent 
$(\alpha/2)$-stable subordinator $S^{(\alpha/2)}$ first  and then killing it upon exiting $D$.  
If we reverse the order of the two operations, that is, we first kill the Brownian motion $W$
upon exiting $D$ and then subordinate the killed Brownian motion $W^D$ using
the independent $(\alpha/2)$-stable subordinator $S^{(\alpha/2)}$, we get the process $Y^{D, (\alpha)}_t:=W^D_{S^{(\alpha/2)}_t}$,
which is called a \textit{subordinate killed Brownian motion}. 
The generator of $X^{(\alpha), D}$ is $-(-\Delta)^{\alpha/2}|_D$, the fractional Laplacian with zero
exterior condition, while the generator of $Y^{D, (\alpha)}$ is $-(-\Delta |_D)^{\alpha/2}$, the fractional
power of the Dirichlet Laplacian.
Subordinate killed Brownian motions are very natural and useful processes. 
For example, it was used in \cite{CS}  as a tool to obtain two-sided estimates for the eigenvalues of the generator of
$X^{(\alpha), D}$. 
The potential theory of subordinate killed Brownian motions has been studied intensively, see
\cite{KSV} and the references therein.
In the PDE literature, the operator $-(-\Delta |_{D})^{\alpha/2}$
also goes under the name of spectral fractional Laplacian, see \cite{BSV} and the references therein. This operator has been of interest to quite a few people in the PDE circle.

The purpose of this paper is to study the small time asymptotic behavior of the spectral
heat content $\widetilde{Q}^{(\alpha)}_D(t)$ with respect to $Y^{D, (\alpha)}$ defined by
$$
\widetilde{Q}^{(\alpha)}_D(t):=\int_D\P_x(Y^{D, (\alpha)}_t\in D)dx.
$$

The main results of this paper are Theorems \ref{thm:SKBM} and \ref{thm:3rd term} below. 
When dealing with stable processes, the notation $\E$ will stand for the expectation of the
process starting from the origin.
An open set $D$ in $\R^{d}$ is said to be a (uniform) $C^{1,1}$ open set
if there are (localization radius) $R_{0}>0$
and $\Lambda_{0}$ such that for every $z\in\partial D$ there exist
a $C^{1,1}$ function $\psi=\psi_{z}:\R^{d}\rightarrow \R$ satisfying $\psi(0, \cdots, 0)=0$,
$\nabla\psi(0)=(0,\cdots,0)$, $|\nabla\psi(x)-\nabla\psi(y)|\leq \Lambda_{0}|x-y|$, and an orthonormal coordinate system
$CS_{z}:y=(y_{1},\cdots,y_{d-1},y_{d}):=(\tilde{y},y_{d})$ with its origin at $z$
such that $B(z,R_{0})\cap D =\{y= (\tilde{y},y_{d})\in B(0,R_{0}) \text{ in } CS_{z} : y_{d}>\psi(\tilde{y})\}$.
In this paper we will call the pair $(R_{0},\Lambda_{0})$ the characteristics of the $C^{1,1}$ open set $D$.

\begin{thm}\label{thm:SKBM}
Let $D=(a,b)$ with $b-a<\infty$ when $d=1$ or $D$ be a bounded $C^{1,1}$ domain when $d\geq 2$. 
Then
$$
\lim_{t\rightarrow 0}\frac{|D|-\widetilde{Q}^{(\alpha)}_{D}(t)}{f_{\alpha}(t)}=
\begin{cases}
\frac{2}{\pi}\Gamma(1-\frac{1}{\alpha})|\partial D|=\frac{2|\partial D|}{\sqrt{\pi}}\E[\left(S_{1}^{(\alpha/2)}\right)^{1/2}] &\mbox{if } 1<\alpha<2,\\
\frac{2}{\pi}|\partial D| &\mbox{if } \alpha=1,\\
\int_{0}^{\infty}(|D|-Q_{D}^{(2)}(u))\frac{\alpha}{2\Gamma(1-\frac{\alpha}{2})}u^{-1-\frac{\alpha}{2}}du &\mbox{if } 0<\alpha<1,
\end{cases}
$$
where $|D|$ is the Lebesgue measure of $D$, $|\partial D|=2$ when $d=1$, $|\partial D|$ is the surface measure of $\partial D$ when $d\geq 2$, and 
$$
f_{\alpha}(t)=
\begin{cases}
t^{1/\alpha} &\text{if } \alpha\in (1,2),\\
t\ln(\frac{1}{t}) &\text{if } \alpha=1,\\
t &\text{if } \alpha\in (0,1).
\end{cases}
$$
\end{thm}

When $\alpha\in (1,2)$ we identify the third term in the small time expansion of 
$\widetilde{Q}^{(\alpha)}_{D}(t)$.

\begin{thm}\label{thm:3rd term}
Suppose that $d\geq 2$, $\alpha\in(1,2)$, and 
$D$ is a bounded $C^{1,1}$ domain in $\R^d$. 
Then the following limit
$$
\lim_{t\rightarrow 0}\frac{\widetilde{Q}^{(\alpha)}_{D}(t)-(|D|-\frac{2|\partial D|}{\sqrt{\pi}}\E[\left(S_{1}^{(\alpha/2)}\right)^{1/2}]t^{1/\alpha})}{t}
$$
exists and its value is given by 
\begin{align*}
&\int_{0}^{1}\left(Q_{D}^{(2)}(u)-(|D|-\frac{2|\partial D|}{\sqrt{\pi}}u^{\frac12})\right)
\frac{\alpha u^{-1-\frac{\alpha}{2}}}{2\Gamma(1-\frac{\alpha}{2})}du\\
&-\frac{1}{\Gamma(1-\frac{\alpha}{2})}\int_{D}\P_{x}(\tau_{D}^{(2)}\leq 1)dx
+\frac{2|\partial D|\alpha}{\sqrt{\pi}(\alpha-1)\Gamma(1-\frac{\alpha}{2})},
\end{align*}
where $\tau_{D}^{(2)}$ is the first time the Brownian motion $W$ exits $D$.
\end{thm}

\begin{remark}
\rm{
We provide some upper and lower bounds for the constant when $D$ is a ball. 
It follows from \cite[Corollary 1.2 (i)]{Savo} and 
\eqref{eqn:3rd term} below
that the constant in Theorem \ref{thm:3rd term} is nonnegative when $\partial D$ is smooth.
It follows from \cite[Theorem 6.2]{BD} that for a bounded $C^{1,1}$ domain $D\subset \R^{d}$, $d\geq 2$ one has 
\beq\label{eqn:remark1}
\left|Q_{D}^{(2)}(u)-\left(|D|-\frac{2|\partial D|u^{1/2}}{\sqrt{\pi}}\right)\right| \leq \frac{10^{d}|D|u}{R_{0}^{2}} \quad u>0,
\eeq
where $(R_{0},\Lambda_{0})$ is the characteristics of the $C^{1,1}$ domain $D$.
Let $D=B(0,r)$ be a ball with radius $r$. In this case we  have $R_{0}=\frac{r}{2}$.
Hence the constant in Theorem \ref{thm:3rd term} is bounded above by
\begin{eqnarray*}
&&\int_{0}^{1}\frac{10^{d}\frac{\omega_{d}}{d}r^{d}u}{(r/2)^{2}}\frac{\alpha u^{-1-\frac{\alpha}{2}}}{2\Gamma(1-\frac{\alpha}{2})}du 
+\frac{2\alpha\omega_{d}}{\sqrt{\pi}(\alpha-1)\Gamma(1-\frac{\alpha}{2})}\\
&=&\frac{2^{2}10^{d} \alpha \omega_{d}}{d\Gamma(1-\frac{\alpha}{2})(2-\alpha)}r^{d-2}
+\frac{2\alpha\omega_{d}}{\sqrt{\pi}(\alpha-1)\Gamma(1-\frac{\alpha}{2})}.
\end{eqnarray*}
On the other hand, by using the trivial bound $\P_{x}(\tau_{D}^{(2)}\leq 1)\leq 1$ and using \cite[Corollary 1.2 (i)]{Savo} we see that the constant in Theorem \ref{thm:3rd term} is bounded below by
$$
-\frac{\omega_{d}r^{d}/d}{\Gamma(1-\frac{\alpha}{2})} +\frac{2\alpha\omega_{d}}{\sqrt{\pi}(\alpha-1)\Gamma(1-\frac{\alpha}{2})}=\frac{\omega_{d}}{\Gamma(1-\frac{\alpha}{2})}\left(\frac{2\alpha }{\sqrt{\pi}(\alpha-1)}-\frac{r^{d}}{d}\right).
$$
Note that when $r=1$ we observe that 
$$
\frac{2\alpha}{\sqrt{\pi}(\alpha-1)}-\frac{1}{d} > \frac{4}{\sqrt{\pi}}-\frac{1}{d} >0 \quad \text{ for all } \alpha\in (1,2). 
$$
}
\end{remark}

Note that, similar to the case of Brownian motion, the first term in the small time expansion of 
$\widetilde{Q}^{(\alpha)}_{D}(t)$ in Theorem \ref{thm:SKBM} involves the volume of the 
domain $D$ and the second term is related to the perimeter $|\partial D|$ of $D$.
In Proposition \ref{prop:different SHs}, we will show that the second terms in 
the small time expansions of 
$\widetilde{Q}^{(\alpha)}_{D}(t)$ and of $Q^{(\alpha)}_{D}(t)$ are different when $\alpha\in (1,2)$ and $D$ is a bounded open interval in $\R^{1}$.

In the Brownian motion case, the third term in the expansion of $Q^{(2)}_{D}(t)$
involves the mean curvature of $D$.
However, the third term in the small time expansion of 
$\widetilde{Q}^{(\alpha)}_{D}(t)$ in Theorem \ref{thm:3rd term} is given by 
a non-explicit expression.
This is probably unavoidable. See the heuristic explanation after the proof of Theorem
\ref{thm:3rd term}.

The organization of the paper is as follows.  
In Section \ref{preliminaries} we fix our notation and recall some basic facts for later use.  
In Section \ref{section:SBM}, the main results, Theorems \ref{thm:SKBM} and \ref{thm:3rd term}, are proved.

In this paper, we use the convention that $c$, lower case or capital, and with or without subscript,
stands for a constant whose value is not important and may change from one appearance to another.

\section{Preliminaries}\label{preliminaries}

\quad
We first collect some basic facts about stable subordinators.
Recall that, for any $\alpha\in (0, 2)$, 
an $(\alpha/2)$-stable subordinator $S_{t}^{(\alpha/2)}$ is a nondecreasing L\'evy process with $S_{0}^{(\alpha/2)}=0$ and 
\beq\label{eqn:SS Laplace}
\E[e^{-\lam S^{(\alpha/2)}_{t}}]=e^{-t\lam^{\alpha/2}}, \qquad \lambda>0, t\ge 0.
\eeq
It is well-known that 
the characteristic exponent of an $(\alpha/2)$-stable subordinator is given by
\beq\label{eqn:CE2}
\Psi(\theta)=|\theta|^{\frac{\alpha}{2}}\left(\cos\frac{\pi \alpha}{4}-i\sin\frac{\pi\alpha}{4}\text{sgn}\theta\right).
\eeq
It follows from \eqref{eqn:SS Laplace} that $S^{(\alpha/2)}_{t}$ and $t^{2/\alpha}
S^{(\alpha/2)}_{1}$ have the same distribution.
The subordinator $S^{(\alpha/2)}$ has a continuous transition density $g^{(\alpha/2)}(t,x)$.
It follows from \cite[(18)]{BKKK} that $g^{(\alpha/2)}(1,x)$ is given by
\beq\label{eqn:SS:HK1}
g^{(\alpha/2)}(1,x)=\frac{1}{\pi}\sum_{n=1}^{\infty}(-1)^{n+1}\frac{\Gamma(1+\frac{\alpha n}{2})}{n!}\sin(\frac{\pi\alpha n}{2})x^{-\frac{\alpha n}{2}-1}, \quad x>0.
\eeq
It follows from Stirling's formula 
$$
\Gamma(1+z)\sim \sqrt{2\pi z} \left(\frac{z}{e}\right)^{z}, \quad z\to \infty,
$$
that for any $\eps>0$ there exists $N$ such that for all $n\geq N$ we have
\beq\label{eqn:Gamma}
\frac{\Gamma(1+\frac{\alpha n}{2})}{n!} \leq \frac{1+\eps}{1-\eps}\frac{\sqrt{2\pi \frac{\alpha n}{2}} \left(\frac{\frac{\alpha n}{2}}{e}\right)^{\frac{\alpha n}{2}}}{\sqrt{2\pi n} \left(\frac{n}{e}\right)^{n}}=\frac{1+\eps}{1-\eps}
\left(\frac{\alpha}{2}\right)^{\frac{\alpha}{2}(n+1)}\left(\frac{n}{e}\right)^{-(1-\frac{\alpha}{2})n},
\eeq
Thus
$$
\sum_{n=N}^{\infty}\left|(-1)^{n+1}\frac{\Gamma(1+\frac{\alpha n}{2})}{n!}\sin(\frac{\pi\alpha n}{2})x^{-\frac{\alpha n}{2}-1}\right|
\le 
\sum_{n=N}^{\infty}\frac{1+\eps}{1-\eps}
\left(\frac{\alpha}{2}\right)^{\frac{\alpha}{2}(n+1)}\left(\frac{n}{e}\right)^{-(1-\frac{\alpha}{2})n}<\infty.
$$
Hence the infinite series in \eqref{eqn:SS:HK1} converges absolutely for all $x>0$. 
Using \eqref{eqn:SS:HK1}, Euler's reflection formula
$$
\Gamma(z)\Gamma(1-z)=\frac{\pi}{\sin(\pi z)}, \quad z\notin \mathbb{Z},
$$
and the absolute convergence of the series in \eqref{eqn:SS:HK1}, we get
\beq\label{eqn:limit2}
\lim_{x\rightarrow\infty}g^{(\alpha/2)}(1,x)x^{1+\frac{\alpha}{2}}=\frac{\Gamma(1+\frac{\alpha}{2})\sin(\frac{\pi\alpha}{2})}{\pi}=\frac{\alpha}{2\Gamma(1-\frac{\alpha}{2})}.
\eeq
(We note in passing that the constant given in \cite[(5.18)]{BBKRSV} is incorrect, due to some typos in transcribing the formula from \cite{Sk}.)
By the scaling property, the transition density $g^{(\alpha/2)}(t,x)$ is equal to $t^{-2/\alpha}g^{(\alpha/2)}(1,\frac{x}{t^{2/\alpha}})$. 
It follows from \eqref{eqn:CE2} and the inverse Fourier transform that for all $x>0$,
\begin{align}\label{eqn:HK ub1}
g^{(\alpha/2)}(1,x)&=(2\pi)^{-1/2}\int_{\R}e^{i\theta x}e^{-|\theta|^{\alpha/2}\left(\cos\frac{\pi \alpha}{4}-i\sin\frac{\pi\alpha}{4}\text{sgn}\theta\right)}d\theta\nn\\
&\leq(2\pi)^{-1/2}\int_{\R}\left|e^{i\theta x}e^{-|\theta|^{\alpha/2}\left(\cos\frac{\pi \alpha}{4}-i\sin\frac{\pi\alpha}{4}\text{sgn}\theta\right)}\right|d\theta\nn\\
&\leq(2\pi)^{-1/2}\int_{\R}e^{-|\theta|^{\alpha/2}\cos\frac{\pi\alpha}{4} }d\theta<\infty. 
\end{align}
On the other hand when $x\geq 1$ it follows from \eqref{eqn:SS:HK1} and \eqref{eqn:Gamma} that
\begin{align}\label{eqn:HK ub2}
g^{(\alpha/2)}(1,x)&=x^{-\frac{\alpha}{2}-1}\frac{1}{\pi}\sum_{n=1}^{\infty}(-1)^{n+1}\frac{\Gamma(1+\frac{\alpha n}{2})}{n!}\sin(\frac{\alpha\pi n}{2})x^{-\frac{\alpha(n-1)}{2}}\nn\\
&\leq x^{-1-\frac{\alpha}{2}}\frac{1}{\pi}\sum_{n=1}^{\infty}\frac{\Gamma(1+\frac{\alpha n}{2})}{n!}x^{-\frac{\alpha(n-1)}{2}}\nn\\
&\leq x^{-1-\frac{\alpha}{2}}\frac{1}{\pi}\sum_{n=1}^{\infty}\frac{\Gamma(1+\frac{\alpha n}{2})}{n!}\leq c_{1}x^{-1-\frac{\alpha}{2}}.
\end{align}
Hence it follows from \eqref{eqn:HK ub1} and \eqref{eqn:HK ub2} that there exists a constant $c_{2}>0$ such that
\beq\label{eqn:SS:HK2}
g^{(\alpha/2)}(t,x) \leq c_{2}\left( t^{-2/\alpha}\wedge \frac{t}{x^{1+\frac{\alpha}{2}}}\right), \quad x>0.
\eeq
We remark here that in case of symmetric stable processes, the transition density also has a 
matching lower bound of the form \eqref{eqn:SS:HK2}. 
However, in the case of stable subordinators, a lower bound similar to the right-hand side of \eqref{eqn:SS:HK2} does not hold (see \cite[Lemma 1]{H}).

The following fact is from \cite[Proposition 2.1]{V1} and will be used in the next section.
\beq\label{eqn:power of subordinator}
\E[(S_{1}^{(\alpha/2)})^{\gamma}]=\frac{\Gamma(1-\frac{2\gamma}{\alpha})}{\Gamma(1-\gamma)}, \quad -\infty<\gamma<\frac{\alpha}{2}.
\eeq

Now we proceed to define for $D$ the killed subordinate Brownian motion and subordinate killed Brownian motion with respect to an isotropic $\alpha$-stable process and their respective spectral heat contents.
Let $W$ be a Brownian motion in $\R^d$ with generator $\Delta$ and 
$S^{(\alpha/2)}_{t}$ an $(\alpha/2)$-stable subordinator independent of $W$. 
Then the subordinate Brownian motion $X^{(\alpha)}$ defined by
$$
X^{(\alpha)}_t:=W_{S^{(\alpha/2)}_{t}}, \quad t\ge 0,
$$
is an isotropic $\alpha$-stable process.
For any domain $D\subset \R^d$, the process $X^{(\alpha), D}$ defined by
$$
X^{(\alpha), D}_t:=
\begin{cases}
X^{(\alpha)}_{t} &\mbox{if } t<\tau^{(\alpha)}_{D},\\
\partial &\mbox{if } t\geq \tau^{(\alpha)}_{D},
\end{cases} \quad t\ge 0,
$$
where 
$$
\tau^{(\alpha)}_{D}:=\inf\{t> 0 : X_{t}^{(\alpha)}\notin D\}, \quad \alpha\in (0,2],
$$
and $\partial$ is a point not contained in $D$ (the cemetery point), is called a killed subordinate Brownian motion,
or more precisely, a killed isotropic $\alpha$-stable process in $D$.  
When $\alpha=2$, $X_{t}^{(2)}$ will be a Brownian motion that will be simply denoted by $W$.
The spectral heat content
of $D$ with respect to $X^{(\alpha)}$ is defined to be
$$
Q^{(\alpha)}_{D}(t):=\int_D\P_x(\tau^{(\alpha)}_{D}>t)dx.
$$

Now let $W^D$ be the killed Brownian motion in $D$. 
The subordinate killed Brownian motion $Y^{D, (\alpha)}$ is defined by
$$
Y^{D, (\alpha)}_{t}:=W^{D}_{S^{(\alpha/2)}_{t}}, \quad t\ge 0.
$$
Let $\zeta^{\alpha}$ be the lifetime of $Y^{D, (\alpha)}$, which is the same as the
first time the process $Y^{D, (\alpha)}$ exits $D$.
The spectral heat content of $D$ with respect to $Y^{D, (\alpha)}$ is defined by 
$$
\widetilde{Q}^{(\alpha)}_{D}(t):=\int_D\P_x(Y^{D, (\alpha)}_t\in D)dx=\int_{D}\P_{x}\left(\zeta^{\alpha}>t\right)dx, \quad t>0.
$$
Note that 
$$
\{\zeta^{\alpha}>t\}=\{\tau^{(2)}_{D}> S^{(\alpha/2)}_{t}\}, \quad t>0.
$$
Hence 
$$
\widetilde{Q}_{D}^{(\alpha)}(t)=\int_{D}\P_{x}\left(\tau_{D}^{(2)}>S^{(\alpha/2)}_{t}\right)dx, \quad t>0.
$$

Note that the following simple relationship is valid
$$
\{\zeta^{\alpha} > t\} =\{\tau_{D}^{(2)}>S_{t}^{(\alpha/2)}\}\subset \{\tau_{D}^{(\alpha)}>t\}, \quad t>0,
$$
which in turn implies
\begin{equation}\label{rel}
\widetilde{Q}_{D}^{(\alpha)}(t)=\int_{D}\P_{x}(\zeta^{\alpha}>t)dx \leq \int_{D}\P_{x}(\tau_{D}^{(\alpha)}>t)dx = Q_{D}^{(\alpha)}(t), \quad t>0.
\end{equation}

We end this section by paraphrasing the explanation given on \cite[page 579]{SV} about
the difference between the processes $X_{t}^{(\alpha), D}$ and $Y_{t}^{(\alpha),D}$. 
Look at a path of the Brownian motion $W$ in $\R^d$, and put a mark on it at all the times given by the subordinator $S_{t}^{(\alpha/2)}$. In this way we observe a trajectory of the process $X_{t}^{(\alpha)}$. The corresponding trajectory of $Y_{t}^{(\alpha),D}$ is given by all the marks on the Brownian path prior to $\tau^{(2)}_{D}$. 
There is a first mark on the Brownian path following the exit time $\tau^{(2)}_{D}$. 
If this mark happens to be in $D$, the process $X_{t}^{(\alpha)}$ has not been killed yet, and the mark corresponds to a point on the trajectory of $X_{t}^{(\alpha),D}$, but not to a point on the trajectory of $Y_{t}^{(\alpha),D}$. 
If, on the other hand, the first mark on the Brownian path following the exit time $\tau^{(2)}_{D}$ happens to be in $D^c$, then trajectories of $Y_{t}^{(\alpha),D}$ and $X_{t}^{(\alpha),D}$ are equal. See the picture on \cite[page 581]{SV} for
an illustration.

\section{Proofs of the main results}\label{section:SBM}

\quad
In this section we prove the main results of this paper, Theorem \ref{thm:SKBM} and Theorem \ref{thm:3rd term}.
First we deal with  the case for $\alpha\in (1,2)$.

\begin{prop}\label{prop:beta12}
Let $\alpha\in (1,2)$. Suppose that $D$ is a bounded open interval when $d=1$ or a bounded  $C^{1,1}$ domain when $d\geq 2$. Then 
$$
\lim_{t\rightarrow0}\frac{|D|-\widetilde{Q}^{(\alpha)}_{D}(t)}{t^{1/\alpha}}=\frac{2|\partial D|}{\sqrt{\pi}}\E\left[\left(S_{1}^{(\alpha/2)}\right)^{1/2}\right].
$$
\end{prop}
\pf
Note that  it follows from \cite[Theorem 1.1]{V2} and \cite[Theorem 6.2]{BD} that
\beq\label{eqn:SH for BM}
\lim_{t\downarrow 0}\frac{|D|-Q_{D}^{(2)}(t)}{t^{1/2}}=\frac{2}{\sqrt{\pi}}|\partial D|.
\eeq
It follows from \eqref{eqn:SH for BM} there exists $\eta>0$ such that
$$
\frac{|D|-Q_{D}^{(2)}(t)}{\sqrt{t}}\leq 1 +\frac{2|\partial D|}{\sqrt{\pi}}, \quad t\in (0,\eta]
$$ 
and if $t\geq \eta$, we have that 
$$
\frac{|D|-Q_{D}^{(2)}(t)}{\sqrt{t}}\leq \frac{|D|}{\sqrt{\eta}}.
$$
Taking $C:=\max\{1+\frac{2|\partial D|}{\sqrt{\pi}}, \frac{|D|}{\sqrt{\eta}}\}$, we get 
\beq\label{eqn:uniform}
\frac{|D|-Q_{D}^{(2)}(t)}{\sqrt{t}}\leq C,\quad  \text{for all } t>0.
\eeq

By the scaling property of $S_{t}^{(\alpha/2)}$ and Fubini's theorem, we have
\begin{align*}
&|D|-\widetilde{Q}_{D}^{(\alpha)}(t)
=\int_{D}\P_{x}(\tau_{D}^{(2)}\leq S_{t}^{(\alpha/2)})dx
=\int_{D}\P_{x}(\tau_{D}^{(2)}\leq t^{2/\alpha}S_{1}^{(\alpha/2)})dx\\
&=\int_{D}\int_{0}^{\infty}\P_{x}(\tau_{D}^{(2)}\leq t^{2/\alpha}u)g^{(\alpha/2)}(1,u)dudx
=\int_{0}^{\infty}\int_{D}\P_{x}(\tau_{D}^{(2)}\leq t^{2/\alpha}u)dxg^{(\alpha/2)}(1,u)du\\
&=\int_{0}^{\infty}\left(|D|-Q_{D}^{(2)}(t^{2/\alpha}u)\right)g^{(\alpha/2)}(1,u)du
=\int_{0}^{\infty}\left(\frac{|D|-Q_{D}^{(2)}(t^{2/\alpha}u)}{t^{1/\alpha}u^{1/2}}\right)t^{1/\alpha}u^{1/2}g^{(\alpha/2)}(1,u)du.
\end{align*}
Hence it follows from \eqref{eqn:SH for BM}, \eqref{eqn:uniform}, \eqref{eqn:power of subordinator} and the Lebesgue dominated convergence theorem that
$$
\lim_{t\rightarrow 0}\frac{|D|-\widetilde{Q}_{D}^{(\alpha)}(t)}{t^{1/\alpha}}=\int_{0}^{\infty}\lim_{t\rightarrow0}\left(\frac{|D|-Q_{D}^{(2)}(t^{2/\alpha}u)}{t^{1/\alpha}u^{1/2}}\right)u^{1/2}g^{(\alpha/2)}(1,u)du=\frac{2|\partial D|}{\sqrt{\pi}}\E\left[\left(S_{1}^{(\alpha/2)}\right)^{1/2}\right].
$$
\qed

Next we deal with the case for $\alpha=1$. We need the following simple lemma. 
\begin{lemma}\label{lemma:beta=1}
For any $\delta>0$, we have
$$
\lim_{t\downarrow 0}\frac{\E\left[\left(S_{1}^{(1/2)}\right)^{\frac{1}{2}}, 0<S_{1}^{(1/2)}<\frac{\delta}{t^{2}} \right]}{\ln(1/t)}=\frac{1}{\sqrt{\pi}}.
$$
\end{lemma}
\pf
It follows from \eqref{eqn:power of subordinator} and an application of Fatou's lemma that
$$
\lim_{t\downarrow 0}\E\left[\left(S_{1}^{(1/2)}\right)^{\frac{1}{2}}, 0<S_{1}^{(1/2)}<\frac{\delta}{t^{2}} \right]=\infty. 
$$
Hence it follows from L'H\^opital's rule, \eqref{eqn:limit2}, 
and the change of variables $x=\delta t^{-2}$ that
\begin{align*}
&\lim_{t\downarrow 0}\frac{\E\left[\left(S_{1}^{(1/2)}\right)^{\frac{1}{2}}, 0<S_{1}^{(1/2)}<\frac{\delta}{t^{2}} \right]}{\ln(1/t)}
=\lim_{t\downarrow 0}\frac{\int_{0}^{\delta/t^{2}}u^{1/2}g^{(1/2)}(1,u)du}{\ln(1/t)}\\
&=\lim_{t\downarrow 0}\frac{(\delta t^{-2})^{1/2}g^{(1/2)}(1,\delta t^{-2})(-2)\delta t^{-3}}{-1/t}
=\lim_{t\downarrow0}2 g^{(1/2)}(1,\delta t^{-2})(\delta t^{-2})^{3/2}\\
&=\lim_{x\uparrow\infty}2 g^{(1/2)}(1,x)x^{3/2}
=\frac{1}{\sqrt{\pi}}.
\end{align*}
\qed

\begin{prop}\label{prop:beta1}
Let $\alpha=1$. Suppose that $D$ is a bounded open interval when $d=1$ or a bounded  
$C^{1,1}$ domain  when $d\geq 2$. Then 
$$
\lim_{t\rightarrow0}\frac{|D|-\widetilde{Q}^{(1)}_{D}(t)}{t\ln(\frac{1}{t})}=\frac{2|\partial D|}{\pi}.
$$
\end{prop}
\pf
As in the proof of Proposition \ref{prop:beta12}. we have
\begin{align}\label{eqn:beta1-1}
&|D|-\widetilde{Q}_{D}^{(1)}(t)
=\int_{0}^{\infty}\left(\frac{|D|-Q_{D}^{(2)}(t^{2}u)}{tu^{1/2}}\right)tu^{1/2}g^{(1/2)}(1,u)du\nn\\
&=\int_{0}^{\delta t^{-2}}\left(\frac{|D|-Q_{D}^{(2)}(t^{2}u)}{tu^{1/2}}\right)tu^{1/2}g^{(1/2)}(1,u)du +\int_{\delta t^{-2}}^{\infty}\left(|D|-Q_{D}^{(2)}(t^{2}u)\right) g^{(1/2)}(1,u)du,
\end{align}
where the value of $\delta$ will be determined later. 

For any $\eps>0$, it follows from \eqref{eqn:SH for BM} that there exists $\delta>0$ such that 
$$
\frac{2|\partial D|}{\sqrt{\pi}} -\eps <\frac{|D|-Q_{D}^{(2)}(t)}{\sqrt{t}} <\frac{2|\partial D|}{\sqrt{\pi}} +\eps, \quad t<\delta.
$$
For this choice of $\delta$ it follows from Lemma \ref{lemma:beta=1} that
\begin{equation}\label{eqn:beta1-2}
\limsup_{t\rightarrow 0}\frac{\int_{0}^{\delta t^{-2}}\left(\frac{|D|-Q_{D}^{(2)}(t^{2}u)}{tu^{1/2}}\right)tu^{1/2}g^{(1/2)}(1,u)du}{t\ln(1/t)}
\leq \left(\frac{2|\partial D|}{\sqrt{\pi}} +\eps\right) \frac{1}{\sqrt{\pi}}.
\end{equation}
Similarly we have
\beq\label{eqn:beta1-3}
\liminf_{t\rightarrow 0}\frac{\int_{0}^{\delta t^{-2}}\left(\frac{|D|-Q_{D}^{(2)}(t^{2}u)}{tu^{1/2}}\right)tu^{1/2}g^{(1/2)}(1,u)du}{t\ln(1/t)}
\geq\left(\frac{2|\partial D|}{\sqrt{\pi}} -\eps\right) \frac{1}{\sqrt{\pi}}.
\eeq

For the second term in \eqref{eqn:beta1-1},  we get from \eqref{eqn:SS:HK2} and the fact $|D|-Q_{D}^{(2)}(t)\leq |D|$ for all $t>0$ that
$$
\int_{\delta t^{-2}}^{\infty}\left(|D|-Q_{D}^{(2)}(t^{2}u)\right) g^{(1/2)}(1,u)du \leq c_{1}|D|\int_{\delta t^{-2}}^{\infty}u^{-3/2}du=c_{2}|D|t
$$
for some constants $c_{1}$ and $c_{2}$. This implies that
\beq\label{eqn:beta1-4}
\limsup_{t\rightarrow 0}\frac{\int_{\delta t^{-2}}^{\infty}\left(|D|-Q_{D}^{(2)}(t^{2}u)\right) g^{(1/2)}(1,u)du}{t\ln(1/t)}=0.
\eeq
Since $\eps$ is arbitrary, the conclusion of the proposition follows from \eqref{eqn:beta1-2}, \eqref{eqn:beta1-3}, and \eqref{eqn:beta1-4}.
\qed

Finally we deal with the case for $\alpha\in (0,1)$. 

\begin{prop}\label{prop:beta01}
Let $\alpha\in (0,1)$. Suppose that $D$ is a bounded open interval when $d=1$ or a bounded  
$C^{1,1}$ domain  when $d\geq 2$. 
Then we have
$$
\lim_{t\downarrow 0}\frac{|D|-\widetilde{Q}_{D}^{(\alpha)}(t)}{t}=\int_{0}^{\infty}\left(|D|-Q_{D}^{(2)}(u)\right)\frac{\alpha}{2\Gamma(1-\frac{\alpha}{2})}u^{-1-\frac{\alpha}{2}}du.
$$
\end{prop}
\pf
Note that by Fubini's theorem,
$$
\frac{1}{t}\left(|D|-\widetilde{Q}^{(\alpha)}_{D}(t)\right)=\frac{1}{t}\int_{D}\P_{x}(\tau_{D}^{(\alpha)}\leq S_{t}^{(\alpha/2)})dx
=\int_{0}^{\infty}\left(|D|-Q_{D}^{(\alpha)}(u)\right)\frac{g^{(\alpha/2)}(t,u)}{t}du.
$$

When $u\geq 1$, it follows from \eqref{eqn:SS:HK2} that
\beq\label{eqn:ub1}
\left(|D|-Q_{D}^{(2)}(u)\right)\frac{g^{(\alpha/2)}(t,u)}{t}\leq c_1|D|u^{-1-\frac{\alpha}{2}}.
\eeq
On the other hand, when $0<u<1$, it follows from \eqref{eqn:uniform} that
$|D|-Q_{D}^{(2)}(u)\leq C u^{\frac12}$. 
Hence from \eqref{eqn:SS:HK2} we have
\beq\label{eqn:ub2}
\left(|D|-Q_{D}^{(2)}(u)\right)\frac{g^{(\alpha/2)}(t,u)}{t}\leq C u^{-\frac{1}{2}-\frac{\alpha}{2}}.
\eeq
Let $\eps>0$ and $\phi_{\eps}\in C_{b}(\R^{1})$ be such that $1_{B(0,\eps)^{c}}\leq \phi_{\eps}\leq 1_{B(0,\frac{\eps}{2})^{c}}$ so that
the function $u\to \left(|D|-Q_{D}^{(2)}(u)\right)\phi_{\eps}(u)\frac{g^{(\alpha/2)}(t,u)}{t}$ is bounded, continuous, and vanishes near zero.
Since $\alpha\in(0,1)$ it follows from \cite[Corollary 8.9]{Sato} and the Lebesgue dominated convergence theorem for any $\eta>0$ there exists $t_{0}>0$ such that
\beq\label{eqn:limitaux}
\left|\int_{0}^{\infty}\left(|D|-Q_{D}^{(2)}(u)\right)\phi_{\eps}(u)\frac{g^{(\alpha/2)}(t,u)}{t}du
-\int_{0}^{\infty}\left(|D|-Q_{D}^{(2)}(u)\right)\phi_{\eps}(u)\frac{\alpha}{2\Gamma(1-\frac{\alpha}{2})}u^{-1-\frac{\alpha}{2}}du\right|<\eta
\eeq
for all $t\leq t_{0}$.
It follows from \eqref{eqn:ub1}, \eqref{eqn:ub2}, and the Lebesgue dominated convergence theorem we have
$$
\lim_{\eps\to 0}\int_{0}^{\infty}\left(|D|-Q_{D}^{(2)}(u)\right)\phi_{\eps}(u)\frac{g^{(\alpha/2)}(t,u)}{t}du=
\int_{0}^{\infty}\left(|D|-Q_{D}^{(2)}(u)\right)\frac{g^{(\alpha/2)}(t,u)}{t}du
$$
uniformly for all $t\leq t_{0}$.
Finally letting $\eps\rightarrow 0$ and using the Lebesgue dominated convergence theorem in \eqref{eqn:limitaux}, we arrive at the conclusion of the proposition.
\qed

\noindent{\bf Proof of Theorem  \ref{thm:SKBM}}.
The proof is an easy consequence of Propositions \ref{prop:beta12}, \ref{prop:beta1}, and \ref{prop:beta01}.
\qed

The second term of the asymptotic expansion of $Q^{(\alpha)}_{D}(t)$ is known when $d=1$
and $D$ is a bouned open interval, see \cite{V2}. We now show that the second terms in
the expansions of  $\widetilde{Q}^{(\alpha)}_{D}(t)$ and of $Q^{(\alpha)}_{D}(t)$ are different when $d=1$ and $\alpha \in (1,2)$. 
\begin{prop}\label{prop:different SHs}
Suppose that $1<\alpha<2$ and $D\subset \R^{1}$ is a bounded open interval. Then 
$$
\lim_{t\rightarrow 0}\frac{|D|-Q^{(\alpha)}_{D}(t)}{t^{1/\alpha}} <\lim_{t\rightarrow 0}\frac{|D|-\widetilde{Q}^{(\alpha)}_{D}(t)}{t^{1/\alpha}}.
$$
\end{prop}
\pf
It is proved in \cite[Theorem 1.1]{V2} that for $\alpha\in (1, 2)$,
$$
\lim_{t\rightarrow 0}\frac{|D|-Q^{(\alpha)}_{D}(t)}{t^{1/\alpha}}=2\E[\overline{X}_{1}^{(\alpha)}].
$$
It follows from \cite[Proposition 2.1]{V2} that
$$
\P(u\leq X^{(\alpha)}_{t})\leq \P(u\leq\overline{X}^{(\alpha)}_{t})\leq 2\P(u\leq X^{(\alpha)}_{t}).
$$
This implies that
\begin{align}\label{eqn:inequality}
&\E[X_{1}^{(\alpha)}, X_{1}^{(\alpha)}>0]=\int_{0}^{\infty}\P(X_{1}^{(\alpha)}\geq u)du
\leq\E[\overline{X_{1}}^{(\alpha)}] =\int_{0}^{\infty}\P(\overline{X}_{1}^{(\alpha)}\geq u)du\nn\\
&\leq \int_{0}^{\infty}2\P(X_{1}^{(\alpha)}\geq u)du=2\E[X_{1}^{(\alpha)}, X_{1}^{(\alpha)}>0].
\end{align}
It is shown in \cite[page 11]{V1} that $\E[X_{1}^{(\alpha)}, X_{1}^{(\alpha)}>0]=\frac{1}{\pi}\Gamma(1-\frac{1}{\alpha})$ and this implies together with Theorem \ref{thm:SKBM}  that
$$
\lim_{t\rightarrow 0}\frac{|D|-Q^{(\alpha)}_{D}(t)}{t^{1/\alpha}}=2\E[\overline{X}_{1}^{(\alpha)}]\leq \frac{4}{\pi}\Gamma(1-\frac{1}{\alpha})=\lim_{t\rightarrow 0}\frac{|D|-\widetilde{Q}^{(\alpha)}_{D}(t)}{t^{1/\alpha}}.
$$

Now we assume  that $\lim_{t\rightarrow 0}\frac{|D|-Q^{(\alpha)}_{D}(t)}{t^{1/\alpha}}= \lim_{t\rightarrow 0}\frac{|D|-\widetilde{Q}^{(\alpha)}_{D}(t)}{t^{1/\alpha}}$.
Then this would imply by \eqref{eqn:inequality} that
$$
2\int_{0}^{\infty}\P(\overline{X}^{(\alpha)}_{1}\geq u)du=2\E[\overline{X_{1}}^{(\alpha)}]\\
=4\E[X_{1}^{(\alpha)}, X_{1}^{(\alpha)}\geq 0]=4\int_{0}^{\infty}\P(X^{(\alpha)}_{1}\geq u)du,
$$
which would imply that $\P(u\leq \overline{X_{1}}^{(\alpha)})=2\P(u\leq X^{(\alpha)}_{1})$ for almost every $u>0$.
But this contradicts \cite[Proposition VIII.4]{Ber} which says  that
$$
\lim_{u\rightarrow \infty}\frac{\P(\overline{X}^{(\alpha)}_{1}\geq u)}{\P(X^{(\alpha)}_{1}\geq u)}=1.
$$
\qed

Combining Theorem \ref{thm:SKBM} with \eqref{rel}, we immediately get the following: when $d\geq 2$, $D$ is a bounded $C^{1,1}$ domain and $\alpha \in (1,2)$, we have
\beq\label{eqn:ub1b}
\limsup_{t\rightarrow 0} \frac{|D|-Q_{D}^{(\alpha)}(t)}{t^{1/\alpha}}\leq \frac{2}{\pi}\Gamma(1-\frac{1}{\alpha})|\partial D|,
\eeq
and when $\alpha=1$
\beq\label{eqn:ub2b}
\limsup_{t\rightarrow 0} \frac{|D|-Q_{D}^{(1)}(t)}{t \ln (\frac{1}{t})}\leq \frac{2}{\pi}|\partial D|.
\eeq
Comparing \eqref{eqn:ub1b} and \eqref{eqn:ub2b} 
with \cite[Theorem 1.3]{V1}, we observe that 
\eqref{eqn:ub1b} and \eqref{eqn:ub2b} are better upper bounds.
We remark that it is conjectured in \cite{V1} that the limits in \eqref{eqn:ub1} 
and \eqref{eqn:ub2} actually exist but this problem is still open. 

Now we establish a three-term small time asymptotic expansion for $\widetilde{Q}^{(\alpha)}_{D}(t)$
when $\alpha\in (1, 2)$. 
First we need the following simple fact. 
\begin{lemma}\label{lemma:auxiliary1}
Let $\alpha\in (1,2)$. Then 
$$
\lim_{t\rightarrow 0}\frac{\E[\left(S^{(\alpha/2)}_{1}\right)^{k/2}, 0< S^{(\alpha/2)}_{1}<t^{-2/\alpha}]}{t^{1-\frac{k}{\alpha}}}
=\frac{\alpha}{(k-\alpha)\Gamma(1-\frac{\alpha}{2})}, \quad k\ge 2.
$$
\end{lemma}
\pf
It follows from \eqref{eqn:power of subordinator} that both the numerator and the denominator diverge to $\infty$ as $t\rightarrow 0$ when $k\geq 2$.
Hence it follows from L'H\^opital's rule and \eqref{eqn:limit2} that
\begin{align*}
&\lim_{t\rightarrow 0}\frac{\E[\left(S^{(\alpha/2)}_{1}\right)^{k/2}, 0< S^{(\alpha/2)}_{1}<t^{-2/\alpha}]}{t^{1-\frac{k}{\alpha}}}
=\lim_{t\rightarrow 0}\frac{\int_{0}^{t^{-2/\alpha}}s^{k/2}g^{(\alpha/2)}(1,s)ds}{t^{1-\frac{k}{\alpha}}}\\
&=\lim_{t\rightarrow 0}\frac{t^{-k/\alpha}g^{(\alpha/2)}(1,t^{-2/\alpha})\frac{-2}{\alpha}t^{-\frac{2}{\alpha}-1}}{(1-\frac{k}{\alpha})t^{-\frac{k}{\alpha}}}
=\lim_{t\rightarrow 0}\frac{2}{k-\alpha}g^{(\alpha/2)}(1,t^{-2/\alpha})t^{-1-\frac{2}{\alpha}}
=\frac{\alpha}{(k-\alpha)\Gamma(1-\frac{\alpha}{2})}.
\end{align*}
\qed

Note that we have
\begin{eqnarray}\label{eqn:3rd term}
&&\widetilde{Q}_{D}^{(\alpha)}(t)-\left(|D|-\frac{2|\partial D|}{\sqrt{\pi}}\E\left[\left(S_{1}^{\alpha/2}\right)^{1/2}\right]t^{1/\alpha}\right)\nn\\
&=&\int_{D}\P_{x}(\tau_{D}^{(2)}\leq S_{t}^{(\alpha/2)})dx-\left(|D|-\frac{2|\partial D|}{\sqrt{\pi}}\E\left[\left(S_{1}^{\alpha/2}\right)^{1/2}\right]t^{1/\alpha}\right)\nn\\
&=&\int_{D}\int_{0}^{\infty}\P_{x}(\tau_{D}^{(2)}\leq st^{2/\alpha})g^{(\alpha/2)}(1,s)dsdx-\left(|D|-\frac{2|\partial D|}{\sqrt{\pi}}\E\left[\left(S_{1}^{\alpha/2}\right)^{1/2}\right]t^{1/\alpha}\right)\nn\\
&=&\int_{0}^{\infty}\left(\int_{D}\P_{x}(\tau_{D}^{(2)}\leq st^{2/\alpha})dx-(|D|-\frac{2|\partial D|}{\sqrt{\pi}}s^{1/2}t^{1/\alpha})\right)g^{(\alpha/2)}(1,s)ds\nn\\
&=&\int_{0}^{\infty}\left(Q_{D}^{(2)}(st^{2/\alpha}) -(|D|-\frac{2|\partial D|}{\sqrt{\pi}}s^{1/2}t^{1/\alpha})\right)g^{(\alpha/2)}(1,s)ds\nn\\
&=&\int_{0}^{t^{-2/\alpha}}\left(Q_{D}^{(2)}(st^{2/\alpha}) -(|D|-\frac{2|\partial D|}{\sqrt{\pi}}s^{1/2}t^{1/\alpha})\right)g^{(\alpha/2)}(1,s)ds\nn\\
&&+\int_{t^{-2/\alpha}}^{\infty}\left(Q_{D}^{(2)}(st^{2/\alpha}) -(|D|-\frac{2|\partial D|}{\sqrt{\pi}}s^{1/2}t^{1/\alpha})\right)g^{(\alpha/2)}(1,s)ds.
\end{eqnarray}

Now we estimate the first expression of \eqref{eqn:3rd term}.
\begin{lemma}\label{lemma:auxiliary2}
Suppose $d\geq 2$ and $\alpha\in (1,2)$. Assume that $D$ is a bounded $C^{1,1}$ domain. Then 
\begin{align*}
&\lim_{t\rightarrow0}\frac{1}{t}\int_{0}^{t^{-2/\alpha}}\left(Q_{D}^{(2)}(st^{2/\alpha}) -(|D|-\frac{2|\partial D|}{\sqrt{\pi}}s^{1/2}t^{1/\alpha})\right)g^{(\alpha/2)}(1,s)ds\\
&=\int_{0}^{1}\left(Q_{D}^{(2)}(u)-(|D|-\frac{2|\partial D|}{\sqrt{\pi}}u^{1/2})\right)\frac{\alpha}{2\Gamma(1-\frac{\alpha}{2})}u^{-1-\frac{\alpha}{2}}du.
\end{align*}
\end{lemma}
\pf
By the change of the variables $u=t^{2/\alpha}s$ and the scaling property of $g^{(\alpha/2)}(t,x)$,
\begin{align}\label{eqn:upper1}
&\frac{1}{t}\int_{0}^{t^{-2/\alpha}}\left(Q_{D}^{(2)}(st^{2/\alpha}) -(|D|-\frac{2|\partial D|}{\sqrt{\pi}}s^{1/2}t^{1/\alpha})\right)g^{(\alpha/2)}(1,s)ds\nn\\
&=\int_{0}^{1}\left(Q_{D}^{(2)}(u)-(|D|-\frac{2|\partial D|}{\sqrt{\pi}}u^{1/2})\right)t^{-\frac{2}{\alpha}-1}g^{(\alpha/2)}(1,t^{-2/\alpha}u)du\nn\\
&=\int_{0}^{1}\left(Q_{D}^{(2)}(u)-(|D|-\frac{2|\partial D|}{\sqrt{\pi}}u^{1/2})\right)\frac{g^{(\alpha/2)}(t,u)}{t}du.
\end{align}
By \eqref{eqn:SS:HK2} we have $\frac{g^{(\alpha/2)}(t,u)}{t}\leq c_{1}u^{-1-\frac{\alpha}{2}}$. 
It follows from \cite[Theorem 6.2]{BD} there exists a constant $c_{2}$ such that 
$$
\left|Q_{D}^{(2)}(u)-(|D|-\frac{2|\partial D|}{\sqrt{\pi}}u^{\frac12})\right|\leq c_{2}u, \quad u>0.
$$
Hence the integrand in \eqref{eqn:upper1} is bounded above by $c_{3}u^{-\frac{\alpha}{2}}$.
Let $\eps>0$ and $\phi_{\eps}\in C_{b}(\R^{1})$ be such that $1_{B(0,\eps)^{c}}\leq \phi_{\eps}\leq 1_{B(0,\frac{\eps}{2})^{c}}$.
Hence it follows from \cite[Corollary 8.9]{Sato} and the Lebesgue dominated convergence theorem that
\begin{align*}
&\lim_{t\rightarrow 0}\int_{0}^{1}\left(Q_{D}^{(2)}(u)-(|D|-\frac{2|\partial D|}{\sqrt{\pi}}u^{1/2})\right)\phi_{\eps}(u)\frac{g^{(\alpha/2)}(t,u)}{t}du\\
&=\int_{0}^{1}\left(Q_{D}^{(2)}(u)-(|D|-\frac{2|\partial D|}{\sqrt{\pi}}u^{1/2})\right)\phi_{\eps}(u)\left(\lim_{t\rightarrow 0}\frac{g^{(\alpha/2)}(t,u)}{t}\right)du\\
&=\int_{0}^{1}\left(Q_{D}^{(2)}(u)-(|D|-\frac{2|\partial D|}{\sqrt{\pi}}u^{1/2})\right)\phi_{\eps}(u)\frac{\alpha}{2\Gamma(1-\frac{\alpha}{2})}u^{-1-\frac{\alpha}{2}}du.
\end{align*}
Letting $\eps\rightarrow 0$, we immediately get the assertion of the lemma by the Lebesgue dominated convergence theorem.
\qed

The two lemmas below are about the second term in \eqref{eqn:3rd term}.

\begin{lemma}\label{lemma:auxiliary3}
Let $\alpha\in (1, 2)$. Then we have
$$
\lim_{t\rightarrow 0}\frac{1}{t}\int_{t^{-2/\alpha}}^{\infty}\left(|D|-Q_{D}^{(2)}(t^{2/\alpha}s)\right)g^{(\alpha/2)}(1,s)ds
=\frac{1}{\Gamma(1-\frac{\alpha}{2})}\int_{D}\P_{x}(\tau_{D}^{(2)}\leq 1)dx.
$$
\end{lemma}
\pf
It follows from Fubini's theorem that
\begin{align*}
&\int_{t^{-2/\alpha}}^{\infty}\left(|D|-Q_{D}^{(2)}(t^{2/\alpha}s)\right)g^{(\alpha/2)}(1,s)ds\\
&=\int_{t^{-2/\alpha}}^{\infty}\int_{D}\P_{x}(\tau_{D}^{(2)}\leq t^{2/\alpha}s)dxg^{(\alpha/2)}(1,s)ds
=\int_{D}\int_{t^{-2/\alpha}}^{\infty}\P_{x}(\tau_{D}^{(2)}\leq t^{2/\alpha}s)g^{(\alpha/2)}(1,s)dsdx.
\end{align*}
It follows from L'H\^opital's rule and \eqref{eqn:limit2} that
\begin{align*}
&\lim_{t\rightarrow 0}\frac{\int_{t^{-2/\alpha}}^{\infty}\P_{x}(\tau_{D}^{(2)}\leq t^{2/\alpha}s)g^{(\alpha/2)}(1,s)ds}{t}\\
&=\lim_{t\rightarrow 0}\P_{x}(\tau_{D}^{(2)}\leq 1)\frac{2}{\alpha}g^{(\alpha/2)}(1,t^{-2/\alpha})t^{-\frac{2}{\alpha}-1}
=\frac{1}{\Gamma(1-\frac{\alpha}{2})}\P_{x}(\tau_{D}^{(2)}\leq 1) .
\end{align*}
Now the result follows from the bounded convergence theorem. 
\qed

\begin{lemma}\label{lemma:auxiliary4}
Let $\alpha\in (1, 2)$. Then 
$$
\lim_{t\rightarrow 0}\frac{1}{t^{1-\frac{1}{\alpha}}}\int_{t^{-2/\alpha}}^{\infty}s^{1/2}g^{(\alpha/2)}(1,s)ds
=\frac{\alpha}{(\alpha-1)\Gamma(1-\frac{\alpha}{2})}.
$$
\end{lemma}
\pf
First note that by \eqref{eqn:power of subordinator} we have
$\lim_{t\rightarrow 0}\int_{t^{-2/\alpha}}^{\infty}s^{1/2}g^{(\alpha/2)}(1,s)ds=0$.
Hence by L'H\^opital's rule and \eqref{eqn:limit2} we have
$$
\lim_{t\rightarrow 0}\frac{1}{t^{1-\frac{1}{\alpha}}}\int_{t^{-2/\alpha}}^{\infty}s^{1/2}g^{(\alpha/2)}(1,s)ds
=\frac{\alpha}{(\alpha-1)\Gamma(1-\frac{\alpha}{2})}.
$$
\qed

\noindent\textbf{Proof of Theorem \ref{thm:3rd term}}
The result follows from Lemmas \ref{lemma:auxiliary1}, \ref{lemma:auxiliary2}, \ref{lemma:auxiliary3}, and \ref{lemma:auxiliary4}.
\qed

Here is a heuristic argument why the third term in the expansion of $\widetilde{Q}^{(\alpha)}_{D}(t)$ involves more than  the mean curvature of $D$.
When $D$ is a bounded smooth domain, the following asymptotic expansion of $Q_{D}^{(2)}(t)$
is well-known (for example see \cite{Savo}):
$$
Q_{D}^{(2)}(t)\sim\sum_{n=0}^{\infty}c_{n}t^{\frac{n}{2}}, \quad \text{as } t\rightarrow 0.
$$
Hence if the  series indeed converges (there exists a case where this series does not converge) and if one could justify the interchange of the integral and the sum, one expects that 
\begin{align*}
&\lim_{t\rightarrow 0}\int_{0}^{1}\left(Q_{D}^{(2)}(u)-(|D|-\frac{2|\partial D|}{\sqrt{\pi}}u^{1/2})\right)\frac{g^{(\alpha/2)}(t,u)}{t}du
=\int_{0}^{1}\left(\sum_{n=2}^{\infty}c_{n}u^{n/2}\right)\lim_{t\rightarrow 0}\frac{g^{(\alpha/2)}(t,u)}{t}du\\
&=\int_{0}^{1}\left(\sum_{n=2}^{\infty}c_{n}u^{n/2}\right)\frac{\alpha}{2\Gamma(1-\frac{\alpha}{2})}u^{-1-\frac{\alpha}{2}}du
=\sum_{n=2}^{\infty}\frac{2c_{n}}{n-\alpha}\frac{\alpha}{2\Gamma(1-\frac{\alpha}{2})}.
\end{align*}
This suggests that even for smooth domains one can not expect that the third term in the expansion of $\widetilde{Q}^{(\alpha)}_{D}(t)$ to involve the mean curvature of $D$ only.
The limit contains the information for all the coefficients of the asymptotic expansion of $Q_{D}^{(2)}(t)$.

\bigskip
\noindent
{\bf Acknowledgment:}
We thank the referee for carefully reading the
manuscript and providing some useful suggestions.

\bigskip
\noindent

\begin{singlespace}

\end{singlespace}

\end{doublespace}

\vskip 0.3truein

{\bf Hyunchul Park}

Department of Mathematics, State University of New York at New Paltz, NY 12561,
USA

E-mail: \texttt{parkh@newpaltz.edu}

\bigskip

{\bf Renming Song}

Department of Mathematics, University of Illinois, Urbana, IL 61801,
USA

E-mail: \texttt{rsong@illinois.edu}

\end{document}